\documentclass[11pt]{article}
\usepackage{amsfonts}

\usepackage{amssymb}
\usepackage{latexsym}
\usepackage{amscd}
\usepackage[centertags]{amsmath}
\usepackage{xypic}
\newtheorem{thm}{Theorem}

\newtheorem{rem}{Remark}

\newcommand{\Enum}{\mathbb{E}}

\newcommand{\Nnum}{\mathbb{N}}

\newcommand{\vep}{\varepsilon}
\newcommand{\qed}{\hfill $\Box$}

\begin{document}
\title{Range-Renewal Speed and
Entropy for I.I.D Models}


\author{Xin-Xing Chen$^{1}$, Jian-Sheng Xie$^{2,3}$
and Jiangang Ying$^{2}$\\
{\footnotesize 1. Department of Mathematics, Shanghai Jiao Tong University,
Shanghai 200240, China}\\
{\footnotesize 2. School of Mathematical Sciences, Fudan University,
Shanghai 200433, China}\\
{\footnotesize 3. Corresponding author. E-mail:
jiansheng.xie@gmail.com}}

\date{}
\maketitle
\begin{abstract}
In this note the relation between the range-renewal speed and entropy for i.i.d. models is
discussed.
\end{abstract}

In \cite{CXY} the authors build an SLLN
\begin{equation}\label{eq: SLLN}
\lim_{n \to \infty} \frac{R_n}{\Enum R_n}=1 \hbox{ almost surely }
\end{equation}
with
\begin{equation}\label{eq: expectation}
\Enum R_n=\sum_x [1-(1-\pi_x)^n]
\end{equation}
for $n$ samples of a discrete distribution $\pi$, where $R_n$ denotes the number of distinct values
of the $n$ samples. In this note we would like to study further the relation between entropy of the
distribution and the range-renewal speed $\Enum R_n$, where the entropy of a (discrete) distribution
$\pi$ is defined as
\begin{equation}
S (\pi) :=\sum_x -\pi_x \cdot \log \pi_x.
\end{equation}
As is already well known, for our i.i.d. model, we always have
$$
\lim_{n \to \infty} \frac{R_n}{n}=0 \hbox{ almost surely.}
$$

But an information of the entropy $S (\pi)$ being finite or infinite
would pose a constriction on the range-renewal speed as the
following:
\begin{thm}\label{thm: entropy-speed}
For our i.i.d. range-renewal model, in general we have
\begin{equation}\label{eq: general}
\lim_{n \to \infty} \frac{R_n}{n}=0
\end{equation}
almost surely. If the entropy $S (\pi)<\infty$, then
\begin{equation}\label{eq: finite-entropy}
\lim_{n \to \infty} \frac{\log n}{n} \cdot R_n=0
\end{equation}
almost surely; Conversely, if the entropy $S (\pi)=\infty$, then
almost surely
\begin{equation}\label{eq: infinite-entropy}
\varlimsup_{n \to \infty} \frac{(\log n)^{1+\vep}}{n} \cdot
R_n=\infty, \quad \forall \vep>0.
\end{equation}
\end{thm}
\noindent \textit{Proof}. \;
We will always assume, for simplicity, that $\pi$ is supported on $\Nnum$ with
$$
\pi_1 \geq \pi_2 \geq \cdots
$$
and we would denote
$$
\varphi^{-1} (n) :=\#\{x: \pi_x>\frac{1}{n}\}
$$
for each $n \geq 1$.

Eq. (\ref{eq: general}) can be proved easily via eq. (\ref{eq: SLLN}) and
(\ref{eq: expectation}) (to prove $\Enum R_n/n \to 0$); it can also be regarded as a
consequence of the main result of \cite{CI78} \cite{CI80} (see also \cite{Athreya85}). Hence
the proof is omitted here.

For (\ref{eq: finite-entropy}), first notice that
\begin{eqnarray*}
\frac{\log n}{n} \cdot \Enum R_n &=& \sum_x [1-(1-\pi_x)^n] \cdot
\frac{\log n}{n}\\
&=& \Bigl\{ \sum_{1/\pi_x \leq n}+\sum_{1/\pi_x>n} \Bigr\} \;
[1-(1-\pi_x)^n] \cdot \frac{\log n}{n} =: I_1+I_2.
\end{eqnarray*}
For the second part of the above equation, we have
\begin{eqnarray*}
I_2 &=& \sum_{1/\pi_x>n} [1-(1-\pi_x)^n] \cdot \frac{\log n}{n}\\
&\leq& \sum_{1/\pi_x>n} \pi_x \cdot \log n \leq \sum_{1/\pi_x>n} -\pi_x \cdot \log \pi_x \to 0.
\end{eqnarray*}
For the first part, choose a large number $N \in \Nnum$ (but still
$N<n$). Noting $\phi (t):=-t \cdot \log t$ is increasing on $(0,
e^{-1})$ (especially on $[1/n, 1/N]$ for all $n>N$), we have
\begin{eqnarray*}
I_1 &=& \sum_{1/\pi_x \leq n} [1-(1-\pi_x)^n] \cdot \frac{\log n}{n} \leq \varphi^{-1} (N) \cdot \frac{\log n}{n}+
\sum_{N<1/\pi_x \leq n} \phi (\frac{1}{n})\\
&\leq& \varphi^{-1} (N) \cdot \frac{\log n}{n}+ \sum_{N<1/\pi_x \leq
n} \phi (\pi_x)\\
&\leq& \varphi^{-1} (N) \cdot \frac{\log n}{n}+ \sum_{x>\varphi^{-1}
(N)} -\pi_x \cdot \log \pi_x.
\end{eqnarray*}
First letting $n \to \infty$ then $N \to \infty$, we get the desired
result.

For (\ref{eq: infinite-entropy}), it's equivalent to $\displaystyle
\varlimsup_{n \to \infty} \frac{(\log n)^{1+\vep}}{n} \cdot \Enum
R_n=\infty$. Suppose on the contrary that there exists some $\vep>0$ such that
$\displaystyle \varlimsup_{n \to \infty} \frac{(\log n)^{1+\vep}}{n} \cdot \Enum
R_n<\infty$. This clearly implies $\displaystyle \varlimsup_{n \to \infty} (\log n)^{1+\vep} \cdot \sum_{\pi_x<1/n}
\pi_x<\infty$ since $\Enum R_n=\sum_x [1-(1-\pi_x)^n]$. We write
$$
a_k:=\#\{x: \frac{1}{k+1}<\pi_x \leq \frac{1}{k}\}, \quad k \geq 1.
$$
Then the above implies $\displaystyle B_n :=(\log n)^{1+\vep} \cdot
\sum_{k=n}^{\infty} \frac{a_k}{k} \leq C$ for some $C>0$ and all $n \geq 1$. Hence
$$
\frac{a_n}{n}=\frac{B_n}{(\log n)^{1+\vep}}-\frac{B_{n+1}}{(\log
(n+1))^{1+\vep}}.
$$
From this we shall derive the following result
\begin{equation}\label{eq: thm-1-1}
\sum_k \frac{a_k}{k} \cdot \log k<\infty,
\end{equation}
which implies $S (\pi)<\infty$, a contradiction. In fact,
\begin{eqnarray*}
\frac{a_n}{n} \cdot \log n &=& \frac{B_n}{(\log n)^{\vep}}
-\frac{B_{n+1}
\cdot \log n}{(\log (n+1))^{1+\vep}}\\
&=& [\frac{B_n}{(\log n)^{\vep}} -\frac{B_{n+1}}{(\log
(n+1))^{\vep}}] +\frac{B_{n+1}
\cdot \log (1+1/n)}{(\log (n+1))^{1+\vep}}\\
&=& [\frac{B_n}{(\log n)^{\vep}} -\frac{B_{n+1}}{(\log
(n+1))^{\vep}}] +O (\frac{1}{n \cdot (\log n)^{1+\vep}}),
\end{eqnarray*}
which surely implies (\ref{eq: thm-1-1}).
\qed

\begin{rem}
\noindent (1) Let
$$
\pi_x :=\frac{C}{x [\log (x+1)]^{\beta+1}}, \quad x =1, 2, \cdots
$$
with $\beta>0$ and $C$ being a normalizing constant. By the results in \cite{CXY} we know
$$
\Enum R_n=O (1) \cdot \frac{n}{(\log n)^{\beta}}
$$
as $n \to +\infty$. When $0<\beta \leq 1$, we always have $S (\pi)=+\infty$, but
$$
\lim_{n \to +\infty} \frac{\log n}{n} \cdot \Enum R_n=\left\{
\begin{array}{rcl}
c \in (0, +\infty), &&\hbox{ if } \beta =1\\
+\infty, &&\hbox{ if } 0<\beta<1
\end{array}
\right.
$$
with $c$ being some positive constant. Therefore the result in (\ref{eq: infinite-entropy})
cannot be strengthened into the one with $\vep=0$;\\

\noindent (2) The $\varlimsup$ in (\ref{eq: infinite-entropy}) cannot be replaced
by $\varliminf$. There exists distributions $\pi$ such that
\begin{equation}\label{eq: anti-infinite-entropy}
S (\pi)=+\infty \hbox{ with } \varliminf_{n \to \infty} \frac{(\log n)^{1+\vep}}{n} \cdot
R_n<\infty, \quad \forall 0<\vep<1.
\end{equation}
\end{rem}
For the part (2) of the above remark, for example, let for any $k\ge 1$,
$$
b_k :=2^{2^k}, \; S_0 :=0, \; S_k :=\sum_{\ell=1}^k\frac{2^{b_{\ell}}}{b_{\ell}}.
$$
And for any $S_{k-1}< x\le S_{k}$, we set $\pi_x :=A \cdot 2^{-b_k}$,
where $A$ is the normalizing constant. Obviously $2<A<4$. It is
easily to see that $S(\pi)=\infty$ since
\begin{eqnarray*}
S(\pi) &=& \sum_{k=1}^\infty\sum_{x=S_{k-1}+1}^{ S_k} \pi_x \log
(\pi_x^{-1}) = \sum_{k=1}^\infty\frac{2^{b_k}}{b_k} \cdot (A \cdot 2^{-b_k})
\log (\frac{2^{b_k}}{A}) \\
&=& \sum_{k=1}^\infty \left[ A \cdot \log 2 -\frac{A}{b_k} \log
A\right]=\infty.
\end{eqnarray*}
The proof of (\ref{eq: anti-infinite-entropy}) is as the following.
For each $k\ge 1$,  let $n_k=2^{2 b_k}$. Then
$A \cdot 2^{-b_{k+1}} < \frac{1}{n_k} \leq A \cdot 2^{-b_k}$ for sufficiently large $k$ and
$\#\{x: \pi_x \geq \frac{1}{n_k}\}=S_k$. And
\begin{eqnarray*}
\Enum R_{n_k}=&\sum \limits_{\pi_x\geq
{n_k}^{-1}} \left[1-(1-\pi_x)^{n_k} \right]+\sum \limits_{\pi_x< {n_k}^{-1}} \left[1-(1-\pi_x)^{n_k} \right]\\
\leq&\sum \limits_{\pi_x\ge {n_k}^{-1}} 1 +{n_k} \cdot \sum \limits_{\pi_x< {n_k}^{-1}}\pi_x
=S_k+n_k \cdot \sum \limits_{\pi_x \le A \cdot 2^{-b_{k+1}}} \pi_x\\
=& S_k +n_k \cdot \sum \limits_{\ell=k+1}^\infty \frac{2^{b_{\ell}}}{b_{\ell}}
\cdot (A \cdot 2^{-b_{\ell}}) =S_k+n_k \cdot \sum \limits_{\ell=k+1}^\infty\frac{1}{b_{\ell}}\\
\leq & \frac{2^{b_k+3}}{b_k}+\frac{2^{2b_k+1}}{b_{k+1}}
=\frac{2^{b_k+3}}{b_k}+\frac{2^{2b_k+1}}{b_k^2}.
\end{eqnarray*}
Fix $0<\vep<1$. Furthermore,
\begin{eqnarray*}
\frac{(\log_2 n_k)^{1+\vep}}{n_k} \cdot \Enum R_{n_k} \le&
\frac{(2b_k)^{1+\vep}}{2^{2b_k}} \cdot \left( \frac{2^{b_k+3}}{b_k}
+\frac{2^{2b_k+1}}{b_k^2} \right)\\
=&\frac{2^{4+\vep}b_k^\vep}{2^{b_k}}+{2^{2+\vep}b_k^{-1+\vep}}\rightarrow0
\end{eqnarray*}
as $k$ tends to infinity. As a result, (\ref{eq: anti-infinite-entropy}) holds.

\noindent{\sl \textbf{Acknowledgements} \quad} The second author would like to thank Prof.
De-Jun Feng for helpful discussions. This work is in part supported by NSFC (No. 11001173,
No. 11271255 and No. 11271077) and the Laboratory of Mathematics for Nonlinear Science,
Fudan University.








\end{document}